\documentclass[11pt,titlepage,oneside,a4paper]{article}
\usepackage{miniquote}
\usepackage{natbib}
\usepackage[resetlabels]{multibib}
\usepackage[utf8]{inputenc}
\newcites{latex}{References}
\usepackage{amsmath}
\usepackage{leftidx}
\usepackage{amsthm}
\usepackage{todonotes}
\usepackage{latexsym}
\usepackage{amssymb}
\usepackage{bbm}
\usepackage[english]{babel}
\usepackage{array}
\usepackage{rotating}
\usepackage{graphicx, color}
\usepackage{multirow}
\usepackage{pslatex}
\usepackage{lscape}
\usepackage[flushmargin, hang]{footmisc}
\usepackage{fancyhdr}
\usepackage{dsfont}
\usepackage{float}
\usepackage{arydshln}
\usepackage{enumitem}
\usepackage{setspace}
\usepackage[T1]{fontenc}
\usepackage{ae, aecompl}
\usepackage{sectsty}
\usepackage{tikz}
\usepackage{comment}
\usepackage[labelfont= bf,sf,small]{caption}
\usepackage[breaklinks, pdfstartview= FitH, colorlinks=true,
linkcolor=blue, citecolor=blue, urlcolor=black]{hyperref}
\usepackage{marginnote}
\usepackage{nameref}
\usepackage{adjustbox}

\allsectionsfont{\sffamily}
\newcommand{\beq}{\begin{equation}}
\newcommand{\eeq}{\end{equation}}
\newcommand{\beqs}{\begin{equation*}}
\newcommand{\eeqs}{\end{equation*}}
\newcommand{\beqn}{\begin{eqnarray}}
\newcommand{\eeqn}{\end{eqnarray}}
\newcommand{\bs}{\begin{eqnarray*}}
\newcommand{\es}{\end{eqnarray*}}
\newcommand{\bi}{\begin{itemize}}
\newcommand{\ei}{\end{itemize}}
\newcommand{\be}{\begin{enumerate}}
\newcommand{\ee}{\end{enumerate}}
\newcommand{\bd}{\begin{description}}
\newcommand{\ed}{\end{description}}

\long\def\symbolfootnote[#1]#2{\begingroup%
\def\thefootnote{\fnsymbol{footnote}}\footnote[#1]{#2}\endgroup}

\begin{document}
\pagestyle{fancy}
\renewcommand{\headrulewidth}{0pt}
\setlength{\topskip}{0pt}
\setlength{\topmargin}{0pt}
\setlength{\headheight}{0pt}
\setlength{\headsep}{0pt}
\setlength{\footskip}{25pt}
\setlength{\skip\footins}{0,5cm}
\setlength{\extrarowheight}{1mm}
\setlength\parindent{0pt}

\rhead[]{}
\lhead[]{}
\chead[]{}
\rfoot[]{}
\cfoot[- \thepage -]{- \thepage\ -}
\lfoot[]{}

\thispagestyle{empty}
\begin{center}

{\Large\bf \textsf{Do algebraic numbers follow Khinchin's
Law?}\symbolfootnote[1]{The authors are grateful to Joshua Lampert, Ulrich Tamm and Michel Waldschmidt for
helpful comments.}
}\\[4ex]

Philipp Sibbertsen,Timm Lampert, Karsten Müller and Michael Taktikos \\[4ex]


\vspace{0.5cm}

\begin{minipage}[h]{\textwidth}
\rule[0ex]{\textwidth}{0.5pt}

{\bf Abstract}

The coefficients of the regular continued fraction for random numbers are distributed by the Gauss-Kuzmin distribution according to Khinchin's law. Their geometric mean converges to Khinchin's constant and their rational approximation speed is Khinchin's speed. It is an open question whether these theorems also apply to algebraic numbers of degree $>2$. Since they apply to almost all numbers it is, however, commonly inferred that it is most likely that non quadratic algebraic numbers also do so. We argue that this inference is not well grounded. There is strong numerical evidence that Khinchin's speed is too fast. For Khinchin's law and Khinchin's constant the numerical evidence is unclear. We apply the Kullback Leibler Divergence (KLD) to show that the Gauss-Kuzmin distribution does not fit well for algebraic numbers of degree $>2$. Our suggestion to truncate the Gauss-Kuzmin distribution for finite parts fits slightly better but its KLD is still much larger than the KLD of a random number. So, if it converges the convergence is non uniform and each algebraic number has its own bound. We conclude that there is no evidence to apply the theorems that hold for random numbers to algebraic numbers.

\vspace{0.5cm}
\end{minipage}
\end{center}



\vspace{0.8cm}

\newpage
\setcounter{footnote}{0}
\setcounter{page}{1}

\section{Introduction}\label{sec:introduction}

The coefficients of the regular continued fraction of a uniformly distributed random number follow Khinchin's Law. Already Gauss had guessed this and \cite{Kuzmin} as well as \cite{Khinchin}, p. 64, theorem 33 proved it and gave an error bound for the convergence.

It is a common folklore that it is likely that algebraic numbers of degree $>2$ also follow Khinchin's Law, cf., e.g., the quote from the backcover of Lang's \cite{Lang2}:

\vspace{0.25cm}

\begin{quotation}
One general idea is that algebraic numbers will exhibit a behaviour that is the same as almost all numbers in a probabilistic sense, except under very specific structural conditions, namely quadratic numbers. Results for almost all numbers (due to Khinchin) show an interplay between calculus and number theory.
\end{quotation}

\vspace{0.25cm}

\cite{Trotter} made $\chi^2$-tests and concluded that the results are as expected. \cite{Bombieri}, p. 141 even go so far to say:

\vspace{0.25cm}

\begin{quotation}
There is no reason to believe that the continued fraction expansions of nonquadratic algebraic irrationals generally do anything other than to faithfully follow Khinchin's Law as detailed below. Indeed experiment suggests that this is even true for parts, short relative to the length
of the period, of the expansions of quadratic irrationals.
\end{quotation}

\vspace{0.25cm}

In contrast, we argue that the belief that Khinchin's Law carries over to nonquadratic algebraic irrationals is not well grounded. Section \ref{sep} provides an overview of our main reasons for this claim. The following sections substantiate the reasons in more detail.  Section \ref{sec:algnumb} discusses continued fractions of algebraic numbers. Section
\ref{sec:goodfit} measures the goodness of fit for our newly proposed truncated Gauss-Kuzmin distribution in comparison to the standard Gauss-Kuzmin distribution. Section \ref{sec:Kc} shows that each algebraic number seems to have its own special behaviour regarding Khinchin's constant. Section \ref{sec:roleofrandom} discusses conjectures regarding the role of Khinchin's Law, Khinchin's constant and Khinchin's approximation speed for algebraic numbers. Section \ref{sec:concl} draws the conclusion.

\section{Reasons for Believing that Algebraic Numbers do not follow Khinchin's Law}\label{sep}

For the following reasons, we question that Khinchin's Law carries over to algebraic numbers with degree $>2$.

\begin{enumerate}
\item Parts of the period of quadratic irrationals show
that the $\chi^2$ tests are not reliable.

\item The Kuzmin distribution is unbounded, while the coefficients of finite parts of the regular continued fraction of algebraic numbers are bounded due to Liouville's theorem. Furthermore, it is known that the
Liouville bounds can be improved. Each number has its own bound and this should be reflected in the distribution or the error bound of its convergence, if it converges. It can presently not be ruled out
that all coefficients of the regular continued fraction of algebraic numbers are bounded. This is an old, open question.

\item Since this is an open question, an account is needed that can also deal with the case of bounded
coefficients. We suggest to truncate the Gauss-Kuzmin distribution for finite parts of nonquadratic algebraic numbers to deal with this. The results show that their KLD  is much larger than the KLD of a random number.


\item Even if the distribution converges to the Gauss-Kuzmin distribution, the numerical evidence shows that the distribution is non uniform. Without clear-cut calculable error bounds this does not help much for finite parts. The error bounds for algebraic numbers, however, are different from those for the random numbers and each algebraic number has its own error bound.

\item The convergence behaviour of the geometric means of the coefficients of the regular continued fractions is very different for each algebraic number. If it converges the convergence is non uniform and each number has its own specific error bound.
For a uniformly distributed random number, however, it converges quite fast to Khinchin's constant.

\item For a typical random number $r_n = \frac{B(n)}{B(n-1) log(B(n-1))} > 1$ has an infinite number of solutions ($\frac{A(n)}{B(n)}$ are the convergents of the random number). But for algebraic numbers numerical evidence by \cite{Trotter}, p. 117, Table B does not support this. For the given 6 algebraic numbers of the 3rd degree there are only 8 such cases for n = 1000 to 3000 and for $\sqrt[3]{7}$ there is not a single such case. This is worse than for the first 1000 convergents, where $\sqrt[3]{2}$ already has 6 such cases. This is numerical evidence that the distribution does not converge to the Gauss-Kuzmin distribution within the margin of the error bounds of random numbers. If it converges, the rate of convergence is slower.

\item There is strong numerical evidence that for any algebraic number $a$ a constant $K(a)>0$ exists with $b_n < K(a) n$ for all natural numbers $n$ with $b_n$ being the coefficients of the regular continued fraction of $a$. So, the arithmetic means of the $b_n$ of algebraic numbers are most probably bounded (like for quadratic irrationals), while for random numbers they diverge to infinity.

\item If such a constant $K(a)>0$ exists, then it can be proven that Levy's error bound for the convergence for random numbers to the Gauss-Kuzmin distribution is not valid for algebraic numbers.

\end{enumerate}

\section{Continued fractions of algebraic numbers} \label{sec:algnumb}

This section considers the distribution of the coefficients of regular continued fractions of algebraic numbers with degree $>2$. \cite{Trotter} aim to show that these coefficients follow a Gauss-Kuzmin distribution. The distribution is defined as follows. Consider a continued fraction expansion of a random number $x$ uniformly distributed in $(0,1)$:

\[x = \frac{1}{k_1 + \frac{1}{k_2+\ldots}}.\]

Then the following holds asymptotically for the distribution of the coefficients $k_n$:

\[\lim_{n \rightarrow \infty} P(k_n = k) = -\log_2 (1 - \frac{1}{(k+1)^2}).\]

The following error bound with a constant $C>0$ was given by\cite{Levy}. With

\[x_n =\frac{1}{k_n + \frac{1}{k_{n+1}+\ldots}}\]

the probability
$\mid P(x_{n} \leq s) - \log_{2}(1+s)\mid $
converges for all $0 <= s <= 1$ to 0 and Levy proved the error bound

\[\mid P(x_n\leq s) - \log_2(1+s)\mid \leq C \cdot 0.7^n.\]

This distribution is known to be the Gauss-Kuzmin distribution. The
approximation of the Gauss-Kuzmin distribution is an asymptotic result
for uniformly distributed random numbers. \cite{Trotter} claim
that this result carries over to algebraic numbers of degree $>2$.
But \cite{Trotter} does not give error bounds for the convergence.

This is problematic as the Gauss-Kuzmin distribution is unbounded by definition while
the coefficients of finite parts of the continued fraction expansion are bounded by
Liouville's theorem.

Liouville's theorem says that for an algebraic number $a$ of a degree $n>2$
there exists a positive constant $C(a)$ with

\[\mid a - \frac{p}{q}\mid > \frac{C(a)}{q^n}\]

for all natural numbers $p$ and $q$.

Liouville's bound is very rough. In general, each algebraic number of degree $> 2$ has its own bound. Because the Thue Siegel Roth theorem (proven by \cite{Roth}) provides only an upper bound for a measure of the speed of the best rational approximation, it is an open question to what extent it can be improved for specific algebraic numbers with effectively computable bounds. \cite{Voutier} investigates this for $\sqrt [3]{2}$ and mentions the following famous result by \cite{Korobov}

 \[\mid \sqrt [3]{2}- \frac{p}{q}\mid > \frac{1}{q^{2.5}}\]

for all natural numbers $p$ and $q$ with the exception of 1 and 4 for $q$.

Hence, effectively computable bounds exist, which are much sharper than Liouville's. These sharper bounds could be used to truncate the Gauss-Kuzmin distribution for algebraic numbers, logarithms of rational numbers and other numbers for which the method with the hypergeometric series, explained by \cite{Voutier}, applies. It is still an open question, whether all periods have such bounds \cite{Waldschmidt}, p. 437, question 2.  A period is a number that can be expressed as an integral of an algebraic function over an algebraic domain. It is difficult to obtain general results because each period has its very own behaviour. It seems to be especially difficult to obtain an optimal lowest upper bound. For certain periods like $\sqrt [3]{2}$  \cite{Voutier}, $\zeta(2)$  \cite{Zudilin} and $\pi$ \cite{Zeilberger} this has been dealt with many times, always slightly improving the results and in all cases research is still going on.  The history for $\pi$ is especially fascinating. \cite{Mahler} started with the upper bound 42 for the effective irrationality measure (in Korobov's formula this is 2.5 as exponent of $q$) and currently \cite{Zeilberger} hold the record with 7.103205334137.


\cite{Waldschmidt}, p. 437 poses a further problem:

\vspace{0.25cm}

\begin{quotation}
A more ambitious goal would be to ask whether real or complex periods behave,
from the point of view of Diophantine approximation by algebraic numbers, like
almost all real or complex numbers.
\end{quotation}

\vspace{0.25cm}

In addition to study differences between random and algebraic numbers, it is an important, though even more difficult endeavour to investigate differences between random numbers and periods. At least for certain classes of periods, it should be possible to obtain results.

Since the Gauss-Kuzmin distribution does not hold for finite parts of the continued
fractions of algebraic numbers, we suggest to take this into account by a truncated version of the Gauss-Kuzmin distribution. The truncated Gauss-Kuzmin distribution for a finite part is defined by the probability function

\begin{equation} \label{eq:truncK}
P_{truncK}(k_n = k) = \frac{P_K(k)}{1 -\log_2(\frac{2 + maxn}{1 + maxn})}. 
\end{equation}

$P_K(k)$ denotes the probability function of the standard Gauss-Kuzmin distribution for $k \leq maxn$ and 0 for $k > maxn$, where $maxn$ is the maximum of the coefficients of the regular continued fraction of the algebraic number for that finite part.

By truncating the Gauss-Kuzmin distribution, we therefore open up the possibility that the
coefficients of the regular continued fraction representation of algebraic
numbers are bounded in general. For infinite unbounded continued fractions our distribution approaches
the standard Gauss-Kuzmin distribution, which is thus also embedded in our
account.

\section{Measuring the goodness of fit of the truncated Kuzmin
distribution} \label{sec:goodfit}

\cite{Trotter} used $\chi^2$ goodness of fit tests to show
that the coefficients of continued fractions of algebraic numbers follow
a Gauss-Kuzmin distribution. They draw their conclusions from
non-rejections of the $\chi^2$ tests. However, as is well known, the
non-rejection of a null hypothesis is no proof for the hypothesis to be
true since a type II error with unknown error probability occurs. In
addition to this, the $\chi^2$ test is known to be unreliable. It is justified mainly because it is easy to use, particularly, when the access to computer power is limited.

Nevertheless the data from the table from \cite{Trotter}, p. 220 are compatible with our approach. The first 5 columns are as follows.

\begin{table}[h]
\centering
    \begin{tabular}{c|c|c|c|c|}
    Number & $\chi^2$ (n=1000) & P & $\chi^2$ (n=3000) & P   \\ \hline
     $\sqrt[3]{2}$ & 4.61 & 0.13 & 5.59 & 0.22\\
     $\sqrt[3]{3}$ & 8.41 & 0.51 & 10.33 & 0.68\\
     $\sqrt[3]{4}$ & 8.47 & 0.51 & 7.71 & 0.44\\
     $\sqrt[3]{5}$ & 8.07 & 0.47 & 9.48 & 0.61\\
     $\sqrt[3]{7}$ & 10.22 & 0.67 & 13.32 & 0.85\\
     \end{tabular}
      \caption{$\chi^2$ test by Lang and Trotter \label{chi2Test}}
\end{table}

P is the approximate probability that $\chi^2$ for a random sample would not be larger. The hypothesis that random behaviour cannot be rejected is in line with our expectation. But we claim that if it converges, the convergence is non uniform. In 4 of 5 cases, P is further away from 0.5 for n=3000 than for n=1000. So it seems to be very difficult to get precise error bounds, if it converges.

As $\chi^2$ test are not very reliable, we revisit this question by using the Kullback Leibler
Divergence (KLD) instead. We do not merely test for one possible
distribution whether it fits the data  but also compare two
different distributions in how they fit the data.

In the following, $\log$ is used and $2$ is the base throughout. The Kullback
Leibler Divergence for the discrete distributions $P$ and $Q$ is defined by

\begin{equation} \label{eq:KLD}
KLD(P,Q) = \sum_{x} P(x) \log(\frac{P(x)}{Q(x)}).
\end{equation}

We have calculated the KLD for 1000 coefficients for algebraic
numbers and the standard Gauss-Kuzmin distribution. Table \ref{KLD} shows the results for some roots of 2.

\begin{table}[h]
\centering
    \begin{tabular}{c|c}
    Algebraic number & KLD \\ \hline
    $\sqrt [3]{2}$ & 0.0955 \\
    $\sqrt[4]{2}$ & 0.0744 \\
    $\sqrt[5]{2}$ & 0.0905 \\
    $\sqrt[6]{2}$ & 0.1103 \\
    $\sqrt[7]{2}$ & 0.1117 \\
    $\sqrt[8]{2}$ & 0.0931 \\
    \end{tabular} \caption{KLD rounded to 4 digits \label{KLD}}
\end{table}

The values in Table \ref{KLD} can be compared to the KLDs of the regular continued fractions of pseudo random numbers. We choose six pseudo random numbers with 1000 decimal digits and calculated the KLDs rounded to 4 digits: 0.0836, 0.0603, 0.0836, 0.0573, 0.0802, 0.0718.
This, as well as our calculations of the KLDs for further algebraic numbers, indicates that the Gauss-Kuzmin distribution fits much better
to the pseudo random numbers than to the algebraic numbers, which shows that the error bounds of the algebraic numbers are worse, if it converges.

Therefore, we suggest to truncate the Gauss-Kuzmin distribution as
defined in equation (\ref{eq:truncK}) for finite parts of the regular continued fraction
of each algebraic number with a
different bound using the properties of the specific number.
Note that this truncation of the distribution is different from the truncation of the Gauss-Kuzmin law, which is used by \cite{Hensley}.

Table \ref{KLD2} provides our calculations of the KLDs for the truncated distribution for 1000 coefficients of the numbers from Table \ref{KLD}.

\begin{table}[h]
\centering
    \begin{tabular}{c|c}
    Algebraic number & KLD \\ \hline
    $\sqrt [3]{2}$ & 0.0953 \\
    $\sqrt[4]{2}$ & 0.0730 \\
    $\sqrt[5]{2}$ & 0.0884 \\
    $\sqrt[6]{2}$ & 0.1102 \\
    $\sqrt[7]{2}$ & 0.1114 \\
    $\sqrt[8]{2}$ & 0.0920 \\
    \end{tabular} \caption{KLD for the truncated Kuzmin distribution
rounded to 4 digits \label{KLD2}}
\end{table}

For the same 6 pseudo random numbers as before we obtain the KLDs
0.0832, 0.0592, 0.0825, 0.0563, 0.0787, 0.0669. The KLDs of the truncated Gauss-Kuzmin distribution are always lower than those of the Gauss-Kuzmin distribution. Note that even the truncated Gauss-Kuzmin distribution still has a quite large
KLD compared to the random numbers. So, it seems likely to consider each number's own distribution and error bounds, if it converges.

The KLDs of the truncated Gauss-Kuzmin distribution also indicate that the Gauss-Kuzmin distribution itself is slightly too large for coefficients of nonquadratic algebraic numbers because the
KLDs of the truncated Gauss-Kuzmin distribution are always better. This
can also be concluded from the fact that both, the standard and the
truncated Gauss-Kuzmin distribution, are only defined for positive
numbers.  Thus, the probability mass of the standard Gauss-Kuzmin
distribution of values greater than our truncation is
redistributed over the finite interval of our truncated version according
to the Gauss-Kuzmin distribution.

Furthermore, it is not even clear whether the truncated distributions converge to the
Gauss-Kuzmin distribution, because it is still an open question whether
the coefficients are bounded or not. \cite{Adamczewski} discuss this and give more sources.

For these reasons, we suggest to use the truncated Gauss-Kuzmin
distribution. This might also be interesting for $\pi$, $\log(2)$ and
other numbers. This is related to Lang's open conjecture that the
approximation speed $B(n)^2$ $\log(B(n))$ is the best possible for
his classical numbers as defined in \cite{Lang1}, p. 635.
This means that if $A(n) / B(n)$ with $B(n) > 0$ is a rational sequence
converging to the classical number $a$, then
for every $\varepsilon > 0$ the product $\mid \hspace{-0.1cm} a-A(n)/B(n)\hspace{-0.1cm}\mid  B(n)^2 (\log(B(n))^{1+\varepsilon}$ 
is always bounded away from 0 \cite{Lang1}, p. 664.

One may go even further and conjecture that this speed is already too fast
and that $\mid \hspace{-0.1cm} a-A(n)/B(n)\hspace{-0.1cm}\mid B(n)^2 (\log(B(n))$ is always bounded away from 0 for Lang's
classical numbers. For quadratic irrationals, for which $B(n)^2$ is the well-known approximation speed, and for Euler's number \cite{Lang2}, p. 76, theorem 5 this is proven but otherwise it is an open question.

\section{Khinchin's constant} \label{sec:Kc}

Another indication that the spead of the rational approximation of each algebraic number of degree $> 2$ should be investigated separately is Khinchin's constant $K0$. \cite{Bailey}, p. 423 write:

\vspace{0.25cm}

\begin{quotation}
It is remarkable that, even though a random fraction’s limiting geometric mean
exists and furthermore equals the Khinchin constant with probability one, not
a single explicit real number (e.g., a real number cast in terms of fundamental
constants) has been demonstrated to have elements whose geometric mean equals K0.
\end{quotation}

\vspace{0.25cm}

Let us abbreviate the truncated Gauss-Kuzmin distributions at $maxn$ by $GKT$ $(maxn)$.
To test the hypothesis that finite parts of the regular continued fraction of algebraic numbers of degree $> 2$ are distributed with $GKT(maxn)$, we first calculate the value $KC(maxn)$ of the limit of the geometric means of the coefficients of the regular continued fraction of random numbers distributed with $GKT(maxn)$:

\[KC(1) = 1, \; KC(maxn) =  2^{GKT(2)} \cdot \ldots \cdot maxn^{GKT(maxn)} \; \mbox{for} \; maxn \geq 2.\]

This sequence converges for $maxn \rightarrow \infty$ monotonously increasing to
Khinchin's constant $K0$ as expected.

This theoretical result for random numbers can now be compared to numerical calculations of the geometrical mean for algebraic numbers. As Table \ref{geometricmean} shows for our test cases, they do not behave like it is expected for the case of random numbers distributed with a truncated Gauss-Kuzmin distribution, since the geometric mean is larger than $KO$.

\begin{table}[h]
\centering
    \begin{tabular}{c|c}
    Number & Geometric mean \\ \hline
    random number & 2.685  \\
    $\sqrt [3]{3}$ & 2.735 \\
    $\sqrt[4]{3}$ & 2.742 \\
    $\sqrt[5]{3}$ & 2.671 \\
    $\sqrt[6]{3}$ & 2.696 \\
    $\sqrt[7]{3}$ & 2.711 \\
    $\sqrt[8]{3}$ & 2.692 \\
    \end{tabular} \caption{Geometric mean for the first 10000 coefficients \label{geometricmean}}
\end{table}

So, if the geometric means $GM$ converge, then this is non uniform and each number has its own error bound. \cite{Bailey}, p. 425 conjecture the following error bound for their example:

 \[
\mid K0 - GM(n)\mid  < \frac{C}{n^{0.5}} 
\]

with $C>0$.

In the case of $\sqrt[4]{3}$ this leads to $C=5.7$ for $n=10000$. This seems reasonable for algebraic numbers, if their geometric means converge. Numerical evidence supports that the convergence for random numbers is much faster.

Similar results are known from other numbers, cf.  \cite{Shiu}, p. 1315 for $\pi$:

\vspace{0.25cm}

\begin{quotation}
...we wish to point out that even if there is a convergence, the rate has to be very slow. It is
easy to see that, with $n=10000$, the change in value of any single partial quotient will have an effect on the third decimal digit for the value of $K(a,n)$. In
fact we found that $K(\pi,10000)$ differs from $K$ by more than $K(\pi,100)$ does.
\end{quotation}

\vspace{0.25cm}

So the convergence behaviour to Khinchin's constant shows that the truncated distributions do not fit better than the Gauss-Kuzmin distribution for algebraic numbers. 

\section{The Role of Randomness} \label{sec:roleofrandom}

The following theorems can be proven for random numbers by probabilistic means:

\begin{description}
\item[T1:] Khinchin's Law holds for the distribution of the coefficients of the regular continued fraction \cite{Khinchin}, pp. 92f.
\item[T2:] Khinchin's constant is the limit of the geometric means of those coefficients \cite{Khinchin}, p. 93.
\item[T3:]  Khinchin's approximation speed $B(n)^2 \log(B(n))$ holds \cite{Khinchin}, p. 69.  The convergents of the random number are $A(n) / B(n)$.
\end{description}

But what role do these theorems play for nonquadratic algebraic numbers?

We state the following conjectures for nonquadratic algebraic numbers:

\begin{description}
\item[C1:] Khinchin's Law is near the distribution and might be valid but precise error bounds for each number are needed as the convergence is non uniform and slower than for random numbers.
\item[C2:] Khinchin's constant is near the geometric means for large $n$ and might be the limit but precise error bounds are needed as the convergence is much slower than for random numbers.
\item[C3:] Khinchin's approximation speed $B(n)^2 \log(B(n))$ does not hold; it is an unreachable upper bound.
\end{description}

According C1, there is an upper bound for the coefficients for finite parts. And even when the Gauss-Kuzmin distribution is truncated for finite parts, it is still only near the real distribution. Using the KLDs, it might be interesting to investigate further what ``near'' means exactly in this context and to search for precise error bounds for each specific number.

According C2, it seems to be very difficult to be more precise here.

C3 means that if $A(n) / B(n)$ with $B(n) > 0$ is a rational sequence
converging to the algebraic number $a$, then
$\mid a-A(n)/B(n)\mid  B(n)^2 (\log(B(n))$ 
always diverges to infinity.
C3 can even be generalised for periods and Lang's classical numbers. It can be proven for quadratic irrationals. For Euler's number, \cite{Adams} proved the speed $B(n)^2 \frac {\log(B(n))} {\log(\log(B(n)))}$,  see also \cite{Lang2}, p. 74. For nonquadratic algebraic numbers, the question is open. It is even open, whether the coefficients are bounded at all. Our calculation seems to indicate that Khinchin's speed is too fast for algebraic numbers, which agrees with the experimental evidence reported by \cite{Trotter}, p. 116:

\vspace{0.25cm}

\begin{quotation}
The tables [for the first 3000 terms of the continued fractions of the cubic numbers] therefore suggest that the type may in fact not be bigger than a constant times the logarithm, and may even be of an order of magnitude smaller than the logarithm.
\end{quotation}

\vspace{0.25cm}

We suggest the new speed $B(n)^2 \frac {\sqrt{\log(B(n))}} {\log(\log(B(n)))}$ and have tested it together with Khinchin's speed for $\sqrt[3]{2}$ for the first 3000 convergents:

\begin{table}[h]
\centering
    \begin{tabular}{c|c|c}
    Range & Khinchin's speed & New speed \\ \hline
    Smaller than 0,5 & 245 & 589 \\
    Between 0,5 and 5 & 23 & 2315 \\
    Between 5 and 10 & 17 & 69\\
    Between 10 and 50 & 114 & 0\\
    Between 50 and 100 & 127 & 0\\
    Larger than 100 & 2474 & 27\\
    \end{tabular} \caption{Speed measurements \label{Speed}}
\end{table}

So, there is strong numerical evidence that Khinchin's speed is too fast, while the new speed fits much better. The 27 measurements larger than 100 merely indicate that the convergents are relatively bad.

Our result basically agrees with \cite{Trotter}, p. 118, Table I and p. 122 Table III, where the authors also investigate other algebraic numbers and come to the conclusion that good measurements with Khinchin's speed are very rare. They
investigate $r_n = \frac{B(n)}{B(n-1) \log(B(n-1))}$ and observe that the values for $r_n$ basically decrease.

Another indication for C3 is the connection of the speed for algebraic numbers with the solutions of diophantine equations. Experimental evidence shows that the speed is slow, when only small solutions exist, cf. \cite{Smart}, p. 135:

\vspace{0.25cm}

\begin{quotation}
Folklore. If a diophantine equation has only finitely many solutions then
those solutions are small in `height' when compared to the parameters of the
equation.
This folklore is, however, only widely believed because of the large amount
of experimental evidence which now exists to support it.
\end{quotation}

\vspace{0.25cm}

From the strongest possible form of the ABC conjecture it follows that only the speed $B(n)^2$ $Exp(\frac {\sqrt[2]{\log(B(n))}} {\log(\log(B(n)))})$ is reachable by nonquadratic algebraic numbers of degree \cite{Frankenhuysen}, p. 46, formula(1.1). Yet, here all ABC equalities are considered and the prime factors contribute to the quality, too. We conjecture that no proper large solutions to resulting equations exist because each fixed algebraic numbers can only have finitely many large ABC hits. Hence, Roth's theorem can be sharpened further. Van Frankenhuysen observes that there is numerical evidence for $O(\log(b_n)) = O(\log(n))$ in his example. This would imply that C3 is valid for this case. 

The data of \cite{Trotter}, p. 118, Table I and p. 122 Table III indicate that for all algebraic numbers $a$ there exists a constant $K(a)$ with $b_n < K(a) n$ for all natural numbers $n$. From their data the following values for $K(a)$ can be calculated:

\begin{table}[h]
\centering
    \begin{tabular}{c|c|c|c}
    Number & $K(a)$ & $n$ of maximum & $r(n)$ \\ \hline
    $\sqrt[3]{2}$ & 14,833 & 36 & 13,8\\
    $\sqrt[3]{3}$ & 5,143 & 119 & 4,6\\
    $\sqrt[3]{4}$ & 25,737 & 579 & 22,0\\
    $\sqrt[3]{5}$ & 160,632 & 19 & 162,7\\
    \end{tabular} \caption{K measurements until n=1000 \label{KMeasurements}}
\end{table}

\begin{table}[h]
\centering
    \begin{tabular}{c|c|c|c}
    Number &$K(a)$ & $n$ of maximum & $r(n)$ \\ \hline
    $\sqrt[3]{2}$ & 10,694 & 1191 & 5,5\\
    $\sqrt[3]{3}$ & 3,971 & 2407 & 3,3\\
    $\sqrt[3]{4}$ & 3,226 & 1974 & 2,7\\
    $\sqrt[3]{5}$ & 15,807 & 1196 & 13,8\\
    \end{tabular} \caption{K measurements from n=1000 until n=3000 \label{KMeasurements2}}
\end{table}

All $K(a)$ are relatively small, which is numerical evidence for C3. It follows that the distribution cannot converge to the Kuzmin distribution within the error bounds for random numbers. This can be seen as follows. Set in Levy's uniform error bound

\[\mid P(x_n<=s) - \log_2(1+s)\mid <= C \cdot 0.7^n\]

the sequence $s_n = \frac {1}{K(a) n}$. Then $P(x_n<=s_n) = 0$ and $\frac {log_2(1+s_n)}{0.7^n}$ is unbounded. So $C$ does not exist and the convergence behaviour must be different from the random numbers.

Furthermore, all $K(a)$ strongly decrease. This is also evidence that C3 is valid. For random numbers, on the other hand, the inequality $b_n > n \cdot \log(n)$ is valid for infinitely many $n$. Thus, for random numbers the arithmetic means of the coefficients diverges to infinity (see \cite{Khinchin}, p.93f), while for algebraic numbers we conjecture that the arithmetic means are bounded due to the numerical evidence.

It is possible to prove C3 for numbers such as quadratic irrationals, $e$ and other numbers with regular continued fractions of Hurwitz' type. For numbers resulting from differential equations $y' = \frac{P(x,y)}{Q(x,y)}$ (or $y''$ etc.) with rational polynomials $P(x,y)$ and $Q(x,y)$, C3 can be proven for many cases using the solution via power series and C fractions, which are then transformed into continued fractions with positive natural numbers $a_n$ and $b_n$. If all $a_{n}=1$, C3 is valid unconditionally. Yet, this cannot always be achieved by this method. One example is the hypergeometric differential equation, where Gauss gave the continued fractions of the solutions. Especially interesting are the confluent hypergeometric functions (see \cite{Jones}, chapter 6.1.2, p. 205-211). In the resulting continued fractions with positive natural $a_n$ and $b_n$, the $b_n$ are linear in $n$ and so Khinchin's speed is too fast. One good example is tanh(z) with $a_0 = z$, $a_n = z^2$, $b_n = 2n + 1$ (formula 6.1.56 of \cite{Jones}). For all natural numbers $z=m$ this is the regular continued fraction and Khinchin's speed is too fast.

Due to the numerical evidence and the fact that there are no known counterexamples to C3, be it as a conjecture applied to algebraic numbers, periods or Lang's classical numbers, we conjecture it for these kind of numbers as a challenge for future research.

The conjectures C1 to C3 compare, in a way, real randomness with pseudo randomness and state that 
the pseudo randomness applying to specific numbers defined by equations is near real randomness but does not reach it exactly. This is known from other generators of pseudo random numbers as well.

\section{Conclusion} \label{sec:concl}

This paper revisits the question whether certain properties of random numbers carry over to algebraic numbers. For random numbers the coefficients of the regular continued fraction asymptotically follow a Gauss-Kuzmin distribution and the convergence has clear error bounds. For algebraic numbers of degree $> 2$ this seems implausible because numerical evidence shows that the convergence is non uniform, if it converges. When the coefficients of the regular continued fractions are bounded, it does not converge at all. We therefore propose a new truncated Gauss-Kuzmin distribution to model the distribution of the coefficients of finite parts of the regular continued fraction of algebraic numbers of degree $>2$. We
apply the Kullback Leibler Divergence to show that our truncated
Gauss-Kuzmin distribution gives a better fit than the standard
Gauss-Kuzmin distribution. This finding is underpinned by simulation
results for a variety of algebraic numbers of degree $>2$. Yet, the KLDs
are still quite large. Furthermore, the convergence behaviour of the
geometric means is very different. It is not even clear, whether the means converge
to Khinchin's constant or not.  Likewise, it is still an open question, whether Khinchin's speed applies to algebraic number of degree $>2$. There is, in fact, strong numerical evidence that Khinchin's speed is too fast. These questions are also interesting to consider for larger classes of numbers like periods or Lang's classical numbers. Probably each such number has its own distribution and error bound. In any case, great care is required when laws that apply to random numbers and, thus, to almost all numbers from a probabilistic point of view, are considered to carry over to specific numbers.








{}


\begin{thebibliography}{}

\bibitem[Adamczewski $\&$ Bugeaud $\&$ Davison(2006)]{Adamczewski} Adamczewski, B., Bugeaud, Y., Davison, L. (2006). Continued Fractions and Transcendental Numbers, \emph{Annales de l'institut Fourier} 56, 2093–2113.

\bibitem[Adams(1966)]{Adams} Adams W. (1966): ``Asymptotic Diophantine approximations to e'', \emph{Proc. Nat. Acad. Sci. U.S.A.} 55 (1966), 28-31. 

\bibitem[Bailey$\&$Borwein$\&$Crandall(1997)]{Bailey} Bailey, D., Borwein, J. and Crandall, R. (1997): ``On the Khinchin Constant'', \emph{Mathematics of Computation} 66(217), 417–431.

\bibitem[Bombieri $\&$ van der Poorten(1975)]{Bombieri} Bombieri, E. and van der Poorten, A. (1975): ``Continued Fractions of Algebraic Numbers'', in: Baker (ed.), \emph{Transcendental Number Theory},
Cambridge University Press, Cambridge, 137-155.

\bibitem[Hensley(1988)]{Hensley} Hensley, D. (1988): ``A Truncated Gauss-Kuzmin Law'', \emph{Transactions of the American Mathematical Society} 306(1), 307-327.

\bibitem[Jones $\&$ Thron(1980)]{Jones} Jones, W.B., Thron, W.J. (1980). Continued fractions. Analytic theory and applications, Encyclopedia of Mathematics and its Applications, 11, Reading, MA: Addison-Wesley.

\bibitem[Khinchin(1963)]{Khinchin} Khinchin, A. (1963): \emph{Continued Fractions}, P.Noordhoff, Groningen.

\bibitem[Korobov(1990)]{Korobov} Korobov, A. (1990): \emph{Continued Fractions and Diophantine Approximations}, Candidate's Dissertation, Moscow State University.

\bibitem[Kuzmin(1928)]{Kuzmin} Kuzmin, R. O.(1928): ``On a Problem of Gauss'', \emph{Doklady Akademii Nauk SSSR}, 375–380.

\bibitem[Lang(1971)]{Lang1} Lang, S. (1971): ``Transcendental Numbers and Diophantine Approximation'',
\emph{Bulletin of the American Mathematical Society} 77(5), 635-677.

\bibitem[Lang $\&$ Trotter(1972)]{Trotter} Lang, S. and Trotter, H. (1972): ``Continued Fractions for some
Algebraic Numbers'', \emph{Journal für die reine und angewandte Mathematik} 255, 112–134; Addendum: ibid. 219–220.

\bibitem[Lang(1995)]{Lang2} Lang, S. (1995): \emph{Introduction to Diophantine Appoximations}, new
expanded edition, New York/Berlin/Heidelberg, Springer.

\bibitem[Lévy(1029]{Levy} Lévy, P. (1929). Sur les lois de probabilité dont dépendant les quotients complets et incomplets d'une fraction continue, \emph{Bulletin de la Société Mathématique de France} 57, 178–194.

\bibitem[Mahler(1953)]{Mahler} Mahler, K. (1953): ``On the Approximation of $\pi$'', \emph{Koninklijke Nederlandse Akademie van Wetenschappen Proceedings} Series A 56(1), 30–42.

\bibitem[Roth(1955)]{Roth} Roth, K. (1955): ``Rational Approximations to Algebraic Numbers and Corrigendum'', \emph{Mathematika} 2, 1-20 and 168.

\bibitem[Shiu(1995)]{Shiu} Shiu, P. (1995), Computation of continued fractions without input values, \emph{Math. Comp.} 64, 1307-1317.

\bibitem[Smart(1998)]{Smart} Smart, N. (1998): \emph{The Algorithmic Resolution of Diophantine Equations}, London Mathematical Society Student Texts 41, London.

\bibitem[Frankenhuysen(1999)]{Frankenhuysen} Van Frankenhuysen, M. (1999): \emph{The ABC conjecture implies Roth's theorem and Mordell's conjecture}, Mat. Contemp. 16 (1999), 45-72.

\bibitem[Voutier(2007)]{Voutier} Voutier, P. (2007): ``Rational Approximation to  $\sqrt [3]{2}$ and other Algebraic Numbers Revisited'', \emph{Journal des Nombres de Bordeaux} 19, 263-288.

\bibitem[Waldschmidt(2006)]{Waldschmidt} Waldschmidt, M. (2006): ``Transcendence of Periods: The State of the Art'', \emph{Pure and Applied Mathematics Quarterly} 2(2), 435—463.

\bibitem[Zeilberger $\&$ Zudilin(2020)]{Zeilberger} Zeilberger, D. and Zudilin, W. (2020): ``The Irrationality Measure of $\pi$ is at most 7.103205334137'', \emph{Moscow Journal of Combinatorics and Number Theory} 9(4), 407–419.

\bibitem[Zudilin(2014)]{Zudilin} Zudilin, W. (2014): ``Two Hypergeometric Tales and a new Irrationality Measure of $\zeta(2)$'', \emph{Annales mathématiques du Québec} 38(1), 101–117.

\end{thebibliography}
\end{document}